\documentclass[12pt]{paper}
\usepackage{amsmath}
\usepackage{hyperref}
\usepackage{amsfonts}
\usepackage{amssymb}
\usepackage{graphics}
\usepackage[pdftex]{graphicx}

\def\bee{\begin{equation}}
\def\eee{\end{equation}}
\def\Li{{\rm {L}i}}
\def\li{{\rm li}}
\begin{document}

\thispagestyle{empty}
\bigskip\bigskip
\centerline{    }
\vskip 2 cm
\centerline{\Large\bf The Skewes number for twin primes:}
\bigskip
\bigskip
\centerline{\Large\bf counting sign changes of $\pi_2(x)-C_2 \Li_2(x)$}

\bigskip\bigskip\bigskip
\centerline{\large\sl Marek Wolf}

\begin{center}
Group of Mathematical Methods in Physics\\
University of Wroc{\l}aw\\
Pl.Maxa Borna 9, PL-50-204 Wroc{\l}aw, Poland\\
e-mail: mwolf@ift.uni.wroc.pl \\
\end{center}

\bigskip\bigskip\bigskip

\begin{center}
{\bf Abstract}\\
\end{center}

\begin{minipage}{12.8cm}
The results of the computer investigation of the sign changes of the
difference between the number of twin primes $\pi_2(x)$ and the
Hardy--Littlewood  conjecture  $C_2\Li_2(x)$ are reported. It turns out that
$d_2(x)=\pi_2(x) - C_2\Li_2(x)$  changes the sign at unexpectedly low values of $x$ and
for $x<2^{48}=2.81\ldots\times10^{14}$ there are 477118 sign changes of this difference.
It is conjectured that the number of sign changes of
$d_2(x)$ for $x\in (1, T)$   is given by $\sqrt T/\log(T)$.
The running logarithmic densities of the sets for which $d_2(x)>0$ and $d_2(x)<0$
are plotted for $x$ up to $2^{48}$.
\end{minipage}

\bigskip\bigskip\bigskip



\bibliographystyle{plain}
\bigskip\bigskip\bigskip

\vfill

\eject

\pagestyle{myheadings}

Let $\pi(x)$ be the number of primes smaller than $x$ and let $\Li(x)$
denote the logarithmic integral:
\begin{equation}
\Li(x)=\int_2^x \frac {du}{ \log(u)}.
\label{Li}
\end{equation}
The Prime Number Theorem tells us that $\Li(x)/\pi(x)$ tends to 1 for
$x\rightarrow\infty$
and the available data (see \cite[Table 14, p. 175]{Ribenboim} or
\cite[Table 5 and 6]{Granville-Races}) show that always
$\Li(x)>\pi(x)$. This last experimental observation was the reason for
the common belief in the past, that the inequality
$\Li(x)>\pi(x)$ is generally valid.
However, in 1914 J.E. Littlewood has shown \cite{Littlewood} (see also \cite{Ellison})
that the difference between the number of primes smaller than $x$ and the
logarithmic integral up to $x$ changes the sign infinitely many times.
The smallest value $x_S$ such that
for the first time $\pi(x_S)\geq \Li(x_S)$ holds is called Skewes number.
We have used ``$\geq$'' to avoid the case of integer value of $\Li(x_S)$, although
we believe that for $n\in \mathbb{N}$ there will be $\Li(n) \notin  \mathbb{N}$,
like we know $\log(n)$ is for $\forall n$ irrational.
In 1933 S. Skewes \cite{SkewesI}
assuming the truth of the Riemann hypothesis  argued that it is certain that
$d(x):=\pi(x)-\Li(x)$ changes sign for some
$x_S < 10^{10^{10^{34}}}$. In 1955
Skewes \cite{SkewesII} has found, without assuming
the Riemann hypotheses, that $d(x)$ changes sign at some
$$
x_S < \exp\exp\exp\exp(7.705)<10^{10^{10^{10^3}}}.
$$
This enormous bound for
$x_S$ was reduced by Cohen and Mayhew \cite{Cohen} to $x_S<
10^{10^{529.7}}$ without using the Riemann hypothesis. In 1966
Lehman \cite{Lehman} has shown that between $1.53\times 10^{1165}$ and
$1.65\times 10^{1165}$ there are more than $10^{500}$ successive integers
$x$ for which $\pi(x)> \Li(x)$. Following the method of Lehman in 1987
H.J.J. te Riele  \cite{Riele} has shown that between $6.62\times 10^{370}$ and
$6.69\times 10^{370}$ there are more than $10^{180}$ successive integers
$x$ for which $d(x)>0$. The lowest present day known estimation of the Skewes number
is around $10^{316}$, see \cite{Bays} and \cite{Demichel2}.

The number of sign changes of the difference $d(x)$ for
$x$ in a given interval $(1, T)$, which is commonly denoted by $\nu(T)$,
see \cite{Ellison}, was discussed for the first time  by A.E. Ingham in 1935
\cite{Ingham1} chapter V, \cite{Ingham2} and next by
S. Knapowski \cite{Knapowski}.
Regarding the number of sign changes of $d(x)$ in the interval
$(1,T)$,  Knapowski \cite{Knapowski} proved that
\begin{equation}
\nu(T)\geq e^{-35}\log \log \log \log T
\label{nu}
\end{equation}
provided $T\geq \exp\exp\exp\exp(35)$. Further results about $\nu(T)$ were obtained
by  J. Pintz \cite{Pintz:III}, \cite{Pintz:IV} and J. Kaczorowski \cite{Kaczorowski:1984a},
\cite{Kaczorowski:1984b}. In particular,  in \cite{Kaczorowski:1984b}  Kaczorowski
proved that there exists such a positive constant $c_3$ that for sufficiently
large $T$ the inequality
\bee
\nu(T) \geq c_3 \log(T)
\eee
holds.  In \cite{Schlage-Puchta-2005} J.-C. Schlage-Puchta proved, assuming
the Riemann Hypothesis, that
\bee
\nu(T)>\frac{\log(T)}{e^{e^{16.7}}}-1.
\eee
More general results on the sign changes can be found in the recent paper
\cite{Kaczorowski_2009}.

In this paper we will look for the analog of the Skewes number for the twin primes,
i.e. pairs of primes separated by 2: \{(3,5), (5,7), (11,13), \ldots, (59, 61), \ldots\}.

Let us denote the number of twin primes pairs $(p, p+2)$ with $p+2<x$
by  $\pi_2(x)$. Then the unproved (see however \cite{Rubinstein})
conjecture B of Hardy and Littlewood \cite{Hardy_and_Littlewood} on
the number of prime pairs $p, p+d$ applied to the case $d=2$ gives, that
\begin{equation}
\pi_2(x) \sim C_2\Li_2(x) \equiv C_2 \int_2^x \frac{u}{\log^2(u)} du,
\label{conj}
\end{equation}
where $C_2$ is called ``twin constant'' and is defined by the following
infinite product:
\begin{equation}
C_2 \equiv 2 \prod_{p > 2} \biggl( 1 - \frac{1}{(p - 1)^2}\biggr) =
1.3203236316937\ldots
\label{stalac2}
\end{equation}

For the first time the conjecture (\ref{conj}) was checked computationally
up to $8\times 10^{10}$ by R. P. Brent
\cite{Brent} who noticed the sign changes of the difference
$\pi_2(x)-C_2 {\rm Li_2}(x)$,  but he did not mention neither the analogy with Skewes
number nor did not count these sign changes.
We analyzed the difference $d_2(x): = \pi_2(x)-C_2\Li_2(x)$ using the
computer  for $x$ up to $T=2^{48}\approx 2.814\times 10^{14}$.
It  took 195 CPU days to reach $T=2^{48}$
on the 64 bits AMD$^{\small{\textregistered}}$ Opteron 2700 MHz processor.

To calculate the integral $\Li_2(x)$ during the main run of the program till $2^{48}$
we have used the 10--point Gauss quadrature \cite{recipes}. This integral was
calculated numerically in successive intervals between
consecutive twins and added to the previous value.  Such a method is not very
time consuming and the number of performed arithmetical operations does not
depend on $x$. There are also power series representations of the logarithmic integral.
We use the following convention for the $\li(x)$ (here {\it v.p. } stands for
French {\it valeur principale} i.e.  Cauchy principal value):
\bee
\li(x)=v.p. \int_0^x \frac{du}{\log(u)} \equiv \lim_{\epsilon \rightarrow 0}\left(
\int_0^{1-\epsilon} \frac{du}{\log(u)} + \int_{1+\epsilon} ^x \frac{du}{\log(u)}\right),
\eee
hence we have $\Li(x)=\li(x)-\li(2)$. Integration by parts gives the asymptotic expansion:
\bee
\li(x) \sim 
\frac{x}{\log( x)}  + \frac{x}{\log^2(x)} + \frac{2x}{\log^3( x)} +
\frac{6x}{\log^4(x)} + \cdots + \frac{n!x}{\log^{n+1}(x)+ \cdots }.
\label{Li_asymptotic}
\eee
which should be cut at $n_0=\lfloor \log(x)\rfloor  $ --- beginning with this index
the following terms are increasing. There is a series giving $\li(x)$ for
all $x>1$ and  quickly convergent which has $n!$ in denominator and $\log^n(x)$
in nominator instead of opposite order in (\ref{Li_asymptotic})
(see \cite[p.126, Entry 14]{BerndtIV})
\bee
\int_\mu^x \frac{du}{\log(u)} = \gamma +
\log \log(x) + \sum_{n=1}^{\infty} {\log^{n}(x)\over n \cdot n!}
\quad {\rm for} ~ x > 1 ~ ,
\label{Li-series}
\eee
where $\gamma=0.5772156649...$ is the Euler-Mascheroni constant and
$\mu=1.451369234883381\ldots$ is the Soldner constant defined by (see
\cite[p.123, eq.(11.3)]{BerndtIV})
\[
\li(\mu)=v.p.\int_0^\mu \frac{du}{\log(u)} = 0.
\]
Even faster converging series was discovered by Ramanujan \cite[p.130, Entry 16]{BerndtIV}:
\bee
\int_\mu^x \frac{du}{\log(u)} = \gamma + \log( \log( x)) + \sqrt{x} \sum_{n=1}^{\infty}
\frac{ (-1)^{n-1} (\log( x))^n} {n! \, 2^{n-1}}
\sum_{k=0}^{\lfloor (n-1)/2 \rfloor} \frac{1}{2k+1} \quad {\rm for} ~ x > 1 ~ .
\label{Li_R}
\eee
Because we have
\[
\Li_2(x)=\Li(x)-\frac{x}{\log(x)}
\]
it is possible to calculate values of $\Li_2(x)$ using the above series. Disadvantage
of these series is that the number of operations (including time consuming
calculation of $\log(x)$) increases with $x$ and is larger
than number of operations needed in the numerical integration.
 	
As for the set of all primes initially the inequality $C_2 {\rm Li}_2(x)> \pi_2(x)$
holds, but it turns out
that there are surprisingly many sign changes of $d_2(x)=\pi_2(x)-C_2\Li_2(x)$
for $x$ in the interval $(1, 2^{48})$. The first sign change of $d_2(x)$ appears
at the twin pair (1369391, 1369393) and up to $T=2^{48}$ there are 477118 sign
changes of $d_2(x)$. We have collected  positions of  all  these sign changes
in one file which is available for downloading from http://www.ift.uni.wroc.pl/\~mwolf/Skewesy\_twins.zip.
Let $\nu_2(T)$ denote, by analogy with usual primes,
the number of sign changes of $d_2(x)$ in the interval
$(1,T)$.   The Table I contains the recorded number of
sign changes of $\pi_2(x) - C_2\Li_2(x)$ up to $T=2^{21}, 2^{22}, \ldots,
2^{48}$.  We have checked the numbers $\nu_2(T)$ up to $T=2^{34}=1.718\times 10^{10}$ independently
calculating the integral $\Li_2(x)$ from the series  (\ref{Li_R}) and these results
are presented in  Table I in the third column and are marked with asterisk.
The first 1274 positions of sign 

\begin{center}
{\sf TABLE {\bf I}}\\
{\small The number of sign changes of $d_2(x)$}\\
\bigskip
\begin{tabular}{||c|c|c|c||c|c|c||} \hline
$T$ & $ \nu_2(T) $ & $ \nu_2(T)(*)$ & $ \sqrt{T}/\log(T)  $ &  $T$ & $ \nu_2(T) $ &  $ \sqrt{T}/\log(T) $\\ \hline
$2^{21}$ &  $ 29 $    &  $ 29  $ &  $         99 $   &     $2^{35}$ &  $ 12682 $ &  $       7641 $ \\ \hline
$2^{22}$ &  $ 29 $    &  $ 29 $ &  $        134 $   &     $2^{36}$ &  $ 23634 $ &  $      10505 $ \\ \hline
$2^{23}$ &  $ 29 $    &  $ 29 $ &  $        182 $   &     $2^{37}$ &  $ 31641 $ &  $      14455 $ \\ \hline
$2^{24}$ &  $ 29 $    &  $ 29 $ &  $        246 $   &     $2^{38}$ &  $ 31641 $ &  $      19905 $ \\ \hline
$2^{25}$ &  $ 29 $    &  $ 29 $ &  $        334 $   &     $2^{39}$ &  $ 31641 $ &  $      27428 $ \\ \hline
$2^{26}$ &  $ 238 $   &  $ 238  $ &  $       455 $  &     $2^{40}$ &  $ 38899 $ &  $      37819 $ \\ \hline
$2^{27}$ &  $ 854 $   &  $ 854  $ &  $       619 $  &     $2^{41}$ &  $ 55106 $ &  $      52180 $ \\ \hline
$2^{28}$ &  $ 1226 $  &  $ 1226  $ &  $      844 $  &     $2^{42}$ &  $ 90355 $ &  $      72037 $ \\ \hline
$2^{29}$ &  $ 1226 $  &  $ 1226 $ &  $       1153 $ &     $2^{43}$ &  $ 161031 $ &  $      99506 $ \\ \hline
$2^{30}$ &  $ 1226 $  &  $ 1226 $ &  $       1576 $ &     $2^{44}$ &  $ 161031 $ &  $     137525 $ \\ \hline
$2^{31}$ &  $ 1226 $  &  $ 1226 $ &  $       2157 $ &     $2^{45}$ &  $ 161031 $ &  $     190168 $ \\ \hline
$2^{32}$ &  $ 2854 $  &  $ 2852  $ &  $       2955 $ &     $2^{46}$ &  $ 405289 $ &  $     263091 $ \\ \hline
$2^{33}$ &  $ 7383 $  &  $ 7381  $ &  $       4052 $ &     $2^{47}$ &  $ 472000 $ &  $     364151 $ \\ \hline
$2^{34}$ &  $ 9115 $  &  $ 9113 $  &  $       5562 $ &     $2^{48}$ &  $ 477118 $ &  $     504258 $ \\ \hline
\end{tabular}
\end{center}
\vskip 0.4cm

\noindent changes of $d_2(x)$ obtained by these two
methods of calculating the integral $\Li_2(x)$ were the same.
The first difference between both methods appears at twin pairs (3067608611, 3067608613)
and  (3067609091, 3067609093), which were  not detected using the more accurate
formula (\ref{Li_R}).  Next twin primes detected by the two methods are the same
until the twin pairs (7809444029,  7809444031). In general, among over 9100
sign changes up to $2^{34}$  there were 17 differences
in the positions of sign changes of $d_2(x)$ obtained by two methods
of calculating the integral $\Li_2(x)$.

The values of $T$ searched by the direct checking are of small magnitude from
the point of view of mathematics, but large for modern computers.

The observed numbers $\nu_2(T)$ behave
somewhat erratically, see Fig.1, in particular there are large gaps without
any change of sign of the $d_2(x)$.  If one assumes the power-like dependence of
$\nu_2(T)$ then the fit by the least square method  gives  the function
$a T^b$, where $a=0.2723\ldots$ and $b=0.4389\ldots$. Instead of such  accidentally
looking parameters of  the pure power-like dependence we suggest the function
$\sqrt T/\log(T)$ as an approximation to $\nu_2(T)$ --- it is a more natural function,
without any free parameters and taking values very close to the least square fit
$a T^b$, see Figure 1. Thus we state the following conjecture:
\begin{equation}
\nu_2(T) \sim \sqrt T/\log(T)~.
\label{conjecture}
\end{equation}
We have picked out function $\sqrt T/\log(T)$ after a few trials and we are not
able to give even heuristic arguments in favour of it.
The conjecture (\ref{conjecture}) is supported by the fact that there are 10 crossings
of the curve $\sqrt T/\log(T)$ with the staircase-like plot of $\nu_2(T)$ obtained
directly from the computer data. The last column  in the Table 1 contains the values
of the function $\sqrt T/\log(T)$. If the conjecture (\ref{conjecture})
is true, then there is infinity of twins.

It seems to be very difficult to gain some analytical insight to why there are so many
sign changes of $\pi_2(x) - C_2\Li_2(x)$.  As (\ref{conj}) is not proved,
hence error term for  it is also not known (for  heuristic approximate formula for
averages of the remainders in the Hardy--Littlewood conjecture B see \cite{KorevaarteRiele}).
The best error term  for Prime Number Theorem
under the Riemann Hypothesis is $|\pi(x)-\Li(x)| = \mathcal{O}(\sqrt{x}\log(x))$.
In the Fig.2 we present the computer data for two functions: the running difference
$d(x)={\rm Li}(x)-\pi(x)$ and the error term:
\bee
\Delta(x) =\max_{2<t<x}\left| \pi(t) - \Li(t) \right|,
\eee
Characteristic oscillations of $d(x)$ are fully described by the explicit formula
for $\pi(x)$, see e.g. \cite[formula (3) and Figure 4]{Granville-Races}.
In the Fig. 3 $|d_2(x)|$ and the error term
\bee
\Delta_2(x) =\max_{2<t<x}\left| \pi_2(t) - C_2\Li_2(t) \right|
\eee
is plotted for $x<2^{48}$. As it is seen from these figures  the behavior of $d(x)$
and $d_2(x)$ is completely different with rapid oscillations of $d_2(x)$ of many orders.
However the functions $\Delta(x)$ and $\Delta_2(x)$
are quite similar: the error term for twins $\Delta_2(x)$ is smaller than $\Delta(x)$
but the difference is not significant:  the power-like fits to $\Delta(x)$
and $\Delta_2(x)$ give:
\bee
\alpha x^\beta,~~~~~~~~\alpha=0.209\ldots,~~\beta=0.45\ldots~~{\rm ~~for~~} \Delta(x)
\eee
\bee
\alpha_2 x^{\beta_2},~~~~~~~~\alpha_2=0.337\ldots,~~\beta_2=0.418\ldots{\rm ~~for~~ } \Delta_2(x).
\eee

Here the  slopes  $\beta \approx \beta_2$ and prefactors
$\alpha$ and $\alpha_2$ are very close.   Thus it seems that the sizes of
the error terms do not account for enormous difference in the value
of Skewes number. In fact all considerations of Skewes, Kaczorowski and others were
based on  existence of explicit formulas and there are no analogs of explicit formulas
for twins.  However  Turan  \cite{Turan}  introduced the following Dirichlet series
with the aim to study twins:
\bee
T(s):= \sum_{n>3}\frac{\Lambda(n-1)\Lambda(n+1)}{n^s} \quad (\Re e\ s > 1),
\label{zeta_twins}
\eee
where $\Lambda(n)$ is the von Mangoldt function:
\bee
\Lambda(n) = \begin{cases}
0 & \mbox{if } n=1 \\
\log p  & \mbox{if }n=p^k \mbox{ for some prime } p \mbox{ and integer } k \ge 1, \\
0 & \mbox{if {\it n} has at least two different prime factors.}
\end{cases}
\eee
In 2004, in a preprint publication \cite{Arenstorf} Arenstorf attempted to prove that
there are infinitely many twins. Arenstorf  tried to continue analytically
$T(s)-C_2/(s-1)$  to  $\Re e\ s=1$,  but shortly after an error in the proof was
pointed out by Tenenbaum \cite{Tenenbaum}.  For  recent progress in the direction of
the proof of the infinite number of twins see \cite{Koreevar}.

 The comparison of Figures 2 and 3 shows,
that $\pi_2(x)\sim C_2\Li_2(x)$ is better than $\pi(x)\sim \Li(x)$ in the sense that
there  are almost half a million points where $d_2(x)$ is zero in the Fig.3 while in the Fig. 2
there are no crossings of $x$ axis at all. This observation can be quantifying
with the notion of the logarithmic density. In \cite{Rubinstein_Sarnak} it was proposed
to  use the logarithmic density to  measure the different biases in the distribution
of prime numbers.  In particular, for the case of the sign changes of
$d(x)$ it was shown that the logarithmic density of the set $\{x:\Li(x)<\pi(x)\}$
defined by
\bee
\delta_{\{x:\Li(x)<\pi(x)\}}=\lim_{x \rightarrow \infty} \frac{1}{\log(x)}\sum_{\substack{2\leq n<x\\ \Li(n)<\pi(n)}}\frac{1}{n}
\label{delta_pi}
\eee
is equal to $\delta_{\{x:\Li(x)<\pi(x)\}}=2.7\ldots\times 10^{-7}$. Hence in some precisely
defined sense the inequality $\Li(x)>\pi(x)$ holds almost everywhere.
Here we will define two logarithmic densities for twin primes as follows:
\bee
\delta_+=\lim_{x \rightarrow \infty} \frac{1}{\log(x)}\sum_{\substack{2\leq n<x\\d_2(n)>0}}\frac{1}{n}
\label{delta1}
\eee

\bee
\delta_-=\lim_{x \rightarrow \infty} \frac{1}{\log(x)}\sum_{\substack{2\leq n<x\\d_2(n)<0}}\frac{1}{n}.
\label{delta2}
\eee
We do not have at our disposal any  formulas like those in \cite{Rubinstein_Sarnak} and
we have to turn  to the brute force numerical calculation of finite size
approximations $\delta_+(x)$ and $\delta_-(x)$ given by expressions
(\ref{delta1}) and (\ref{delta2}) without limit operation $\lim_{x \rightarrow \infty}$.
In these computation we have used positions of all sign changes collected earlier.
The resulting running logarithmic densities are  plotted in Figure 4.
The sum for $\delta_-(x)$ starts from 1/5, because 5 is the end of the first
twin primes pair. It is a reason why the plot of $\delta_-(x)$ in Fig.4 starts
from about 0.67.  Up to $x=2^{31}$  the data for Figure 4 was
obtained by direct summing of the harmonic sums, for $x>2^{31}\approx 2.15\times 10^9$
the incredible accurate approximation   \cite{detemple1993}, \cite[pp. 76-78]{Havil03}:
\bee
\sum_{k=n}^m \frac{1}{k}=\log\left(m+\frac{1}{2}\right)- \log\left(n-\frac{1}{2}\right) + {\mathcal{O}}\left(\frac{1}{n^2}\right)
\eee
was used (the implied in $\mathcal{O}$ constant is much smaller than 1). For
$n\approx 10^9$ the error made by using the above formula is of the
order $10^{-18}$.  To calculate the harmonic
series up to $x=2.8\times 10^{14}$ directly by adding all numbers $1/n$ would take
from one to a few months of CPU time, depending on the processor. The plots presented
in Fig.4 suggest following the conjecture
\bee
\delta_+=\delta_-=\frac{1}{2}.
\eee

\bigskip

The difference of many  hundreds of orders between  values of $x$ such that
$\pi(x)- {\rm Li}(x)$ and  $\pi_2(x)-C_2{\Li_2}(x)$ changes the sign for the  first
time is astonishing. We can give an example from physics. Let us make the mapping:
sign changes of $d(x)$ correspond to energy levels of hydrogen and
sign changes of $d_2(x)$ correspond to the spectrum  of helium. Then
ground states of hydrogen and of helium will correspond to $x_S$ and
first  sign change of $d_2(x)$ accordingly. The experiments show that
the energies  of the ground states of the hydrogen and helium are
-13.6 eV and  -79 eV respectively and do not differ by hundreds of orders!

\bigskip

\noindent{\bf Acknowledgement}:  We would like to thank Prof. A. Jadczyk for  reading
the manuscript and for helpful comments. The author is also very much thankful to the
referees for their valuable remarks and suggestions.

\newpage

\begin{figure}
\includegraphics[width=\textwidth, angle=0, scale=1]{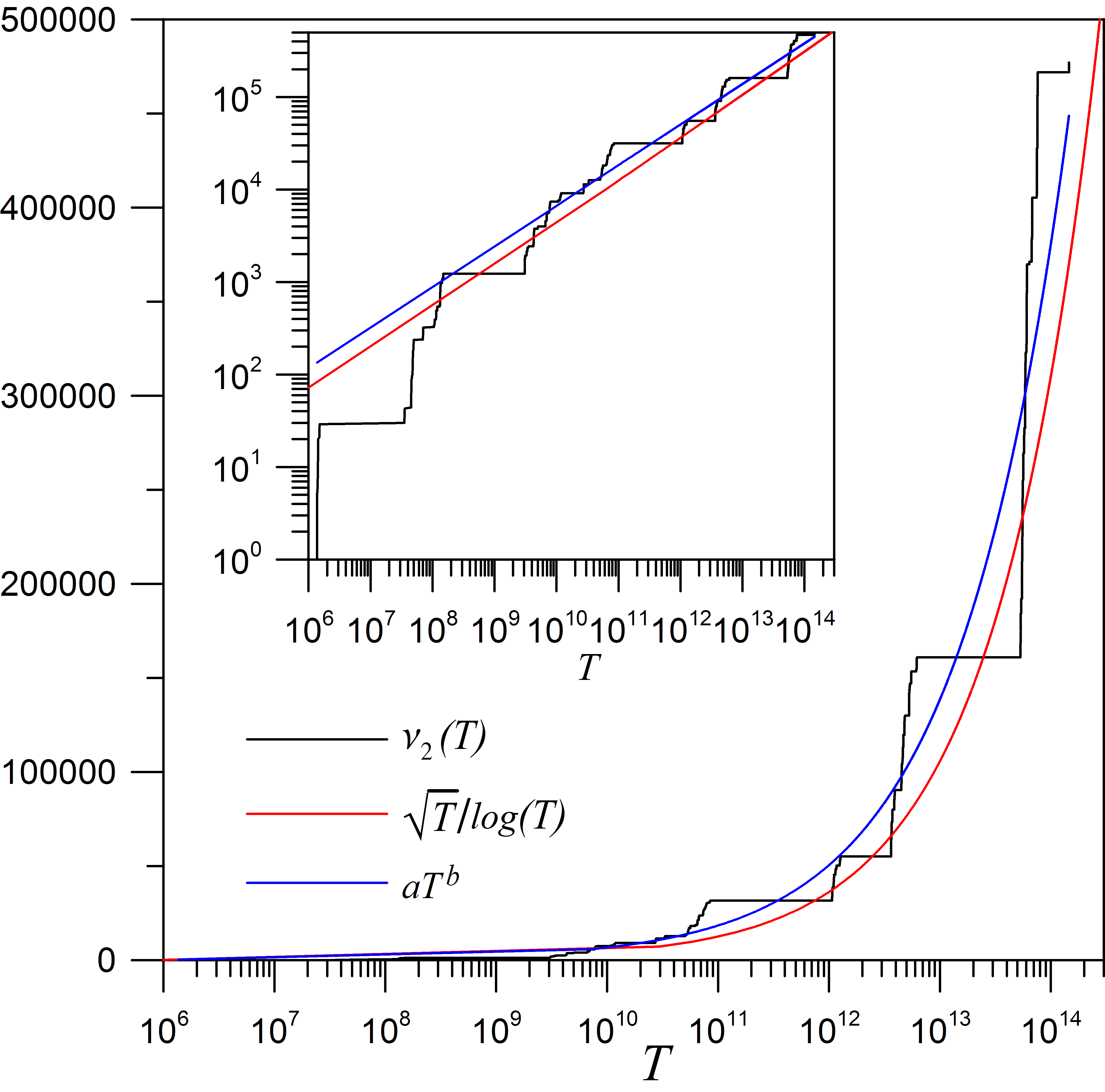} \\
Fig.1 The plot showing the comparison of the actual values of $\nu_2(T)$
found by a computer search with the conjecture  (\ref{conjecture}).
There are 10 crossing of the function $\nu_2(T)$ and $\sqrt{T}/\log(T)$
in this plot up to $2^{48}$.  All 477118 sign changes of $d_2(x)$ are
plotted.  In the inset plot on the double logarithmic scale is presented.\\
\end{figure}

\newpage

\begin{figure}
\includegraphics[width=\textwidth, angle=0, scale=1]{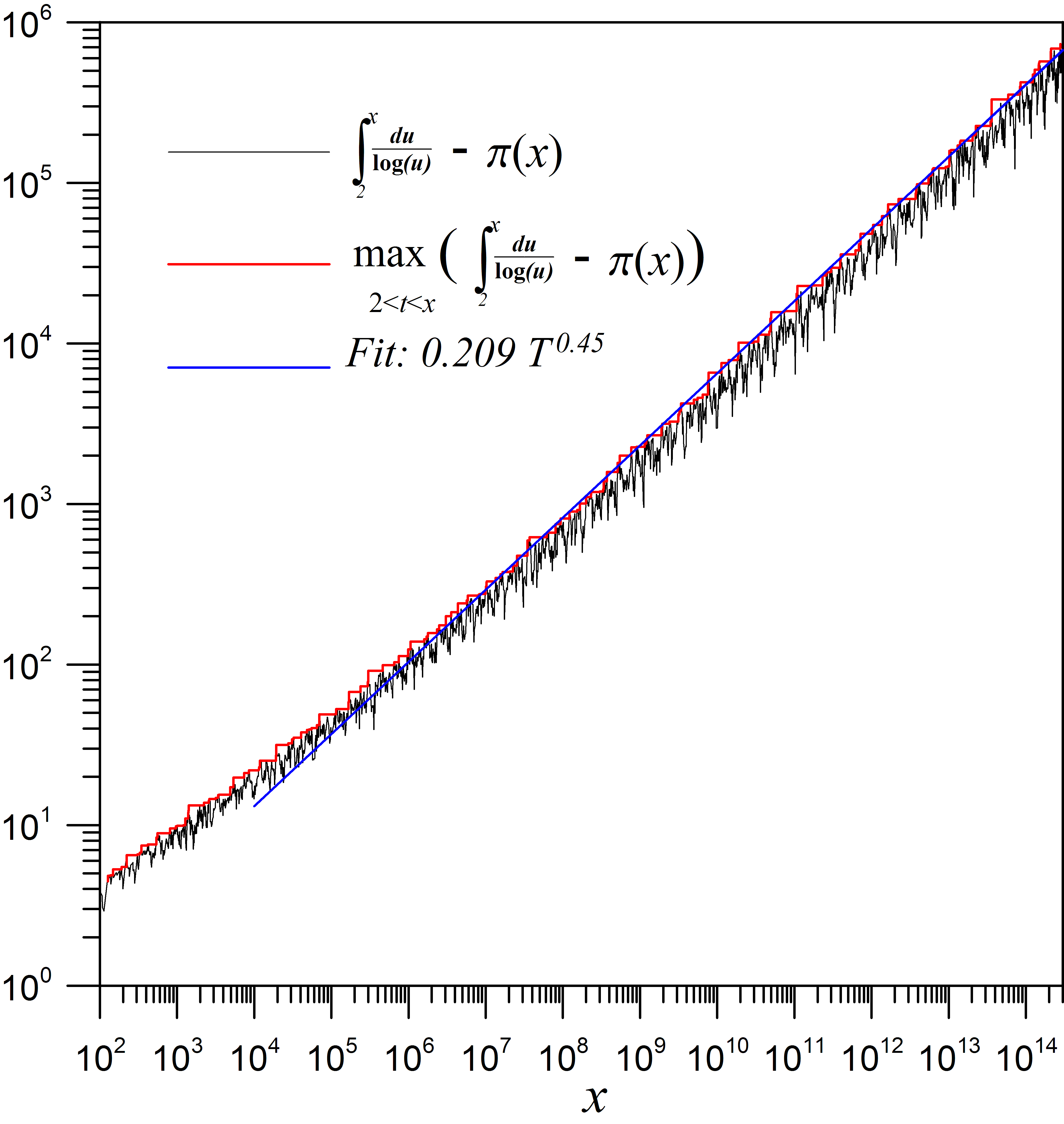} \\
Fig.2 The  plot of $d(x)$ and  error term $\Delta(x)$. The power fit was made
for $10^6<x<2^{48}$. The first crossing of the axis $x$ will appear around $10^{316}$.
\end{figure}

\newpage

\begin{figure}
\includegraphics[width=\textwidth, angle=0, scale=1]{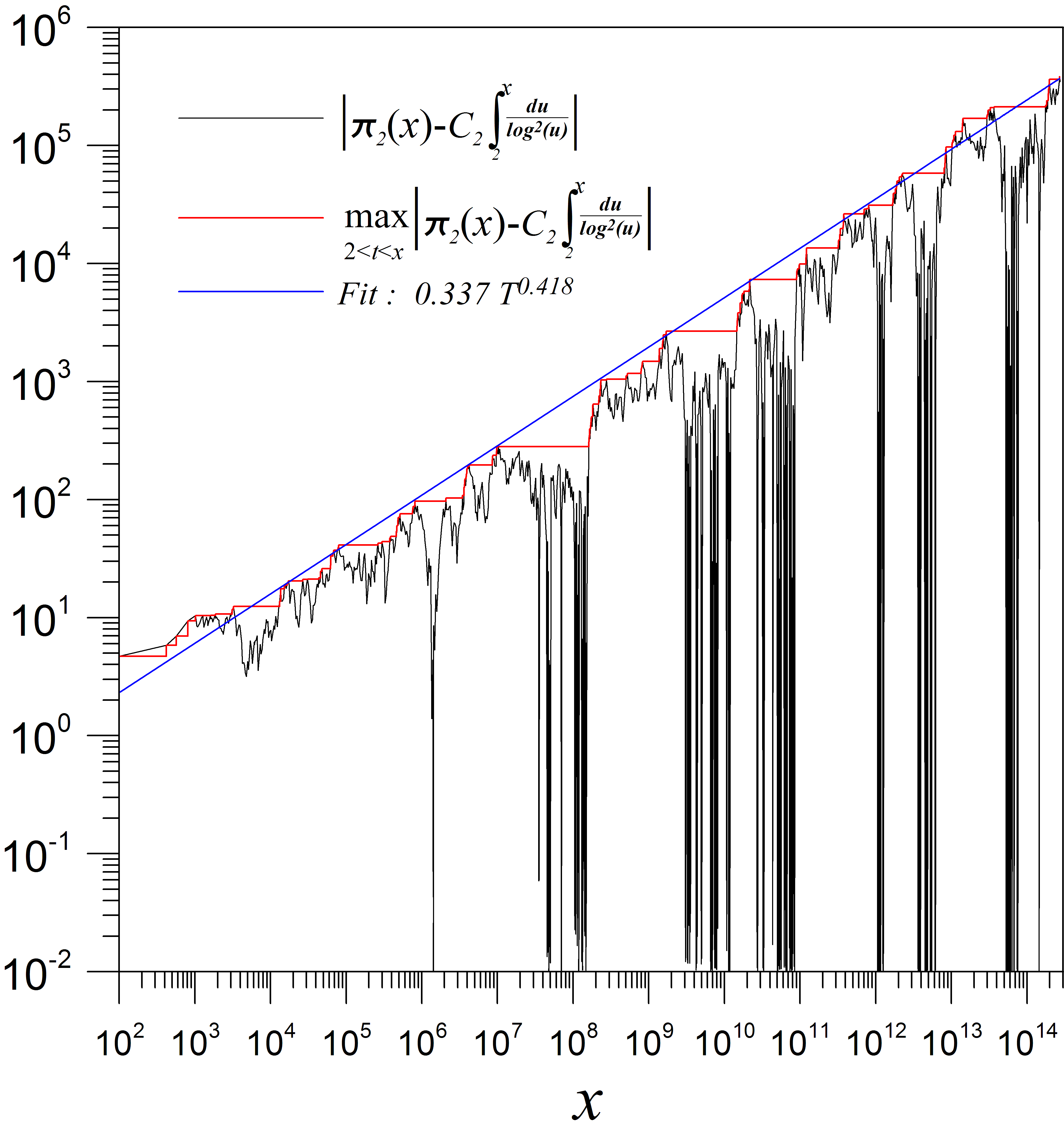}\\
Fig.3 The plot  of $|d_2(x)|$ and error term   $\Delta_2(x)$.
Sign changes of the $d_2(x)$ and  values smaller than $10^{-2}$ were artificially
set to $10^{-2}$. In blue the  power-like fit $0.337 \times x^{0.418}$
to $\Delta_2(x)$ obtained by the least-square method is plotted.
\end{figure}

\newpage

\begin{figure}
\includegraphics[width=\textwidth, angle=0, scale=1]{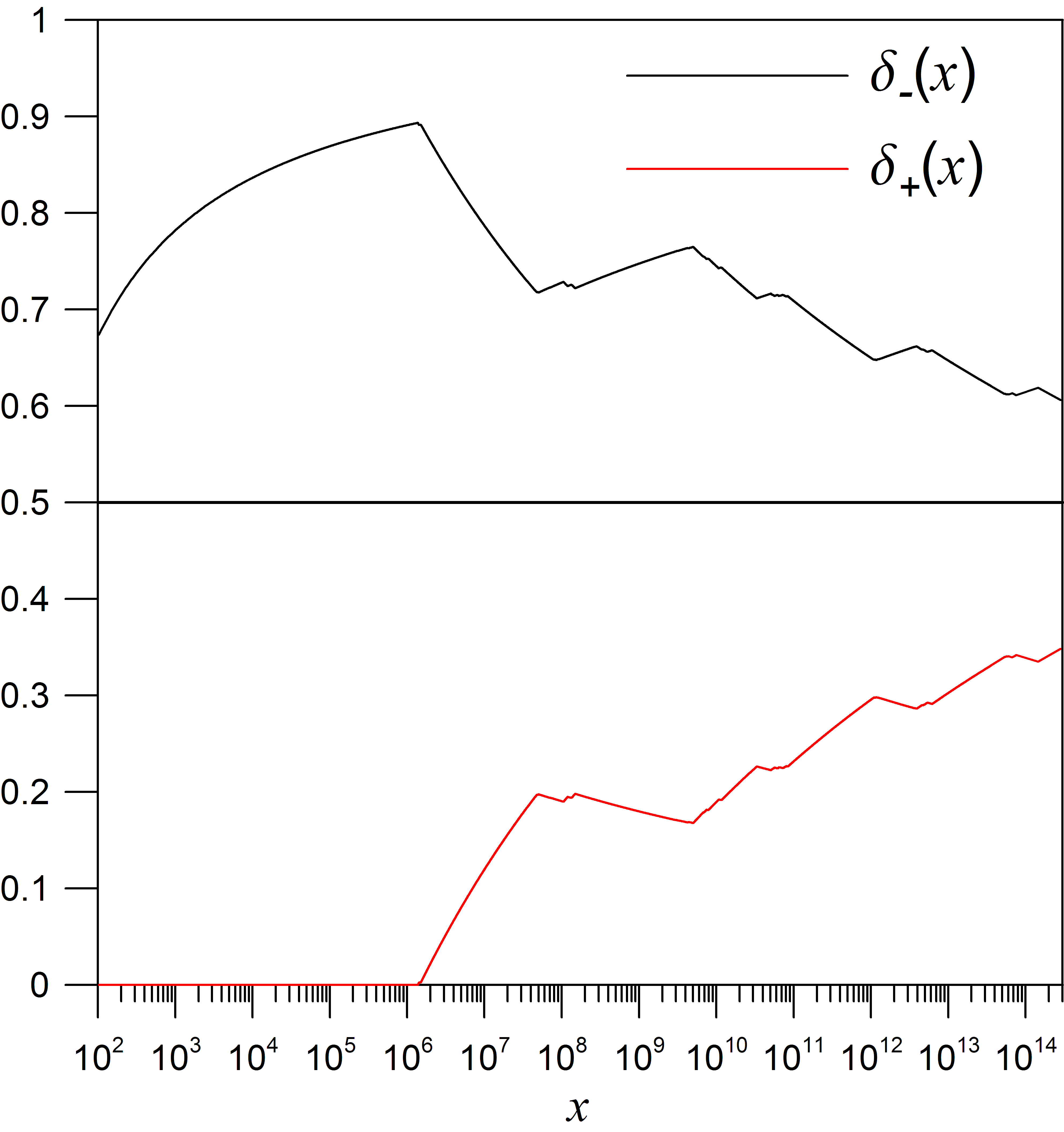}\\
Fig.4  The plots of the running logarithmic densities $\delta_+(x), \delta_-(x)$
defined in the text. Each plot consists of
28025 points: the values of $\delta(x)$'s were  recorded at the progression
$x=100\times(1.001)^n$.
\end{figure}

\end{document}